\documentclass[twocolumn]{autart}    

\usepackage{graphicx}          

\usepackage{amsmath,amsfonts,amssymb,color}      
\usepackage{hyperref}
\hypersetup{
    colorlinks=true,
    linkcolor=blue,
    citecolor=red,
    filecolor=magenta,      
    urlcolor=cyan,
           }                                                                      
\usepackage{epsfig}
\usepackage{psfrag}
\usepackage{algorithm}
\usepackage[noend]{algpseudocode}
\usepackage{array}
\usepackage{epstopdf}
\usepackage{cite}
\usepackage{mathtools}
\usepackage[utf8]{inputenc}
\usepackage{subfig}
\usepackage{caption}
\usepackage{bm}
\usepackage{wrapfig}

\newcommand{\W}{\mathcal{W}}
\newcommand{\B}{\mathcal{B}}
\newcommand{\A}{\mathcal{A}}
\newcommand{\D}{\mathcal{D}}

\newcommand{\BE}{\mathbb{E}}
\newcommand{\R}{\mathbb{R}}
\newcommand{\CS}{\mathcal{S}}
\newcommand{\CF}{\mathcal{F}}

\newcommand{\Z}{\mathbb{Z}}

\newcommand{\CN}{\mathcal{N}}
\newcommand{\bw}{\bm{\mathrm{w}}}
\newtheorem{observation}{Observation}
\DeclarePairedDelimiter\norm{\lVert}{\rVert}
\newcommand{\tr}{\text{Tr}}

\begin{document}

\begin{frontmatter}

\title{Maximization of Nonsubmodular Functions under Multiple Constraints with Applications} 


\author[Hust]{Lintao Ye}\ead{yelintao93@hust.edu.cn},    
\author[Hust]{Zhi-Wei Liu}\ead{zwliu@hust.edu.cn},
\author[Hust]{Ming Chi}\ead{chiming@hust.edu.cn},
\author[Pu]{Vijay Gupta}\ead{gupta869@purdue.edu},               

\address[Hust]{School of Artificial Intelligence and Automation, Huazhong University of Science and Technology, Wuhan 430074, China}  



\address[Pu]{The Elmore Family School of Electrical and Computer Engineering, Purdue University, West Lafayette, IN 47907, USA}             


\begin{abstract}                          
We consider the problem of maximizing a monotone nondecreasing set function under multiple constraints, where the constraints are also characterized by monotone nondecreasing set functions. We propose two greedy algorithms to solve the problem with provable approximation guarantees. The first algorithm exploits the structure of a special class of the general problem instance to obtain a better time complexity. The second algorithm is suitable for the general problem. We characterize the approximation guarantees of the two algorithms, leveraging the notions of submodularity ratio and curvature introduced for set functions. We then discuss particular applications of the general problem formulation to problems that have been considered in the literature.  We validate our theoretical results using numerical examples. 
\end{abstract}

\end{frontmatter}

\section{Introduction}
We study the problem of maximizing a set function over a ground set $\CS$ in the presence of $n$ constraints, where the constraints are also characterized by set functions. Specifically, given monotone nondecreasing set functions\footnote{A set function $f:2^{\CS}\to\R_{\ge0}$ is monotone nondecreasing if $f(\A)\le f(\B)$ for all $\A\subseteq\B\subseteq\CS$.} $f:2^{\CS}\to\R_{\ge0}$ and $h_i:2^{\CS}\to\R_{\ge0}$ for all $i\in[n]\triangleq\{1,\dots,n\}$ with $n\in\Z_{\ge1}$, we consider the following constrained optimization problem:
\begin{equation}\tag{P1}
\label{eqn:obj org}
\begin{split}
&\max_{\A\subseteq \CS}f(\A)\\
&s.t.\ h_i(\A\cap \CS_i)\le H_i,\ \forall i\in[n],
\end{split}
\end{equation}
where $\CS=\cup_{i\in[n]}\CS_i$, with $\CS_i\subseteq \CS$ for all $i\in[n]$, and $H_i\in\R_{\ge0}$ for all $i\in[n]$. By simultaneously allowing nonsubmodular set functions\footnote{A set function $f:2^{\CS}\to\R_{\ge0}$ is submodular if and only if $f(\A\cup\{v\})-f(\A)\ge f(\B\cup\{v\})-f(\B)$ for all $\A\subseteq\B\subseteq\CS$ and for all $v\in\CS\setminus\B$.} (in both the objective and the constraints) and multiple constraints (given by upper bounds on the set functions), \eqref{eqn:obj org} generalizes a number of combinatorial optimization problems (e.g., \cite{khuller1999budgeted,iyer2012algorithms,leskovec2007cost,kulik2009maximizing,bian2017guarantees,das2018approximate}). Instances of \eqref{eqn:obj org} arise in many important applications, including sensor (or measurement) selection for state (or parameter) estimation (e.g., \cite{jawaid2015submodularity,zhang2017sensor,ye2021parameter}), experimental design (e.g., \cite{krause2008near,bian2017guarantees}), and data subset (or client) selection for machine learning (e.g., \cite{das2011submodular,durga2021training,ye2021client}). As an example, the problem of sensor selection for minimizing the error covariance of the Kalman filter (studied, e.g., in \cite{jawaid2015submodularity,ye2020complexity}) can be viewed as a special case of \eqref{eqn:obj org} where the objective function $f(\cdot)$ is defined on the ground set $\CS$ that contains all candidate sensors, and is used to characterize the state estimation performance of the Kalman filter using measurements from an allowed set $\A\subseteq\CS$ of selected sensors. The constraints corresponding to $h_i(\cdot)$ for all $i\in[n]$ represent, e.g., budget, communication, spatial and energy constraints on the set of selected sensors (e.g., \cite{mo2011sensor,ye2020complexity,prasad2022policies}). We apply our results to two  specific problems in Section~\ref{sec:applications}.

In general, \eqref{eqn:obj org} is NP-hard (e.g., \cite{feige1998threshold}), i.e., obtaining an optimal solution to \eqref{eqn:obj org} is computationally expensive. For instances of \eqref{eqn:obj org} with a monotone nondecreasing submodular objective function $f(\cdot)$ and a single constraint, there is a long line of work for showing that greedy algorithms yield constant-factor approximation ratios for \eqref{eqn:obj org} (e.g., \cite{nemhauser1978analysis,khuller1999budgeted,calinescu2011maximizing}). However, many important applications that can be captured by the general problem formulation in \eqref{eqn:obj org} do not feature objective functions that are submodular (see, e.g., \cite{krause2008near,das2011submodular,zhang2017sensor,elenberg2018restricted,ye2020complexity}). For instances of \eqref{eqn:obj org} with a nonsubmodular objective function $f(\cdot)$, it has been shown that the greedy algorithms yield approximation ratios that depend on the problem parameters (e.g., \cite{das2011submodular,bian2017guarantees,Tzoumas2020LQG}). As an example, the approximation ratio of the greedy algorithm provided in~\cite{bian2017guarantees} depends on the submodularity ratio and the curvature of the objective function in \eqref{eqn:obj org}. 

Moreover, most of the existing works consider instances of \eqref{eqn:obj org} with a single constraint on the set of the selected elements, e.g., a cardinality, budget, or a matroid constraint. The objective of this paper is to relax this requirement of a single simple constraint being present on the set of selected elements. For instance, in the Kalman filtering based sensor scheduling (or selection) problem described above, a natural formulation is to impose a separate constraint on the set of sensors selected at different time steps, or to consider multiple constraints such as communication and budget constraints on how many sensors can work together simultaneously. Thus, we consider the problem formulation \eqref{eqn:obj org}, where the objective function is a monotone nondecreasing set function and the constraints are also characterized by monotone nondecreasing set functions. We do not assume that the objective function and the functions in the constraints are necessarily submodular. We propose approximation algorithms to solve \eqref{eqn:obj org}, and provide theoretical approximation guarantees for the proposed algorithms leveraging the notions of curvature and submodularity ratio. 

Our main contributions are summarized as follows. First, we consider instances of \eqref{eqn:obj org} with $\CS_i\cap\CS_j=\emptyset$ for all $i,j\in[n]$ ($i\neq j$), and propose a parallel greedy algorithm with time complexity $O((\max_{i\in[n]}|\CS_i|)^2)$ that runs for each $i\in[n]$ in parallel. We characterize the approximation guarantee of the parallel greedy algorithm, leveraging the submodularity and curvature of the set functions in the instances of \eqref{eqn:obj org}. Next, we consider general instances of \eqref{eqn:obj org} without utilizing the assumption on mutually exclusive sets $S_{i}$. We propose a greedy algorithm with time complexity $O(n|\CS|^2)$, and characterize its approximation guarantee. The approximation guarantee of this algorithm again depends on the submodularity ratio and curvature of the set functions and the solution returned by the algorithm. Third, we specialize these results to some example applications and evaluate these approximation guarantees by bounding the submodularity ratio and curvature of the set functions. Finally, we validate our theoretical results using numerical examples; the results show that the two greedy algorithms yield comparable performances that are reasonably good in practice. A preliminary version of this paper was presented in \cite{ye2021client}, where only the parallel greedy algorithm was studied for a special instance of \eqref{eqn:obj org}.

\textbf{Notation}
For a matrix $P\in\R^{n\times n}$, let $P^{\top}$, $\tr(P)$, $\lambda_1(P)$ and $\lambda_n(P)$ be its transpose, trace, an eigenvalue with the largest magnitude, and an eigenvalue with the smallest magnitude, respectively. A positive definite matrix $P\in\R^{n\times n}$ is denoted by $P\succ 0$. Let $I_n$ denote an $n\times n $ identity matrix. For a vector $x$, let $x_i$ (or $(x)_i$) be the $i$th element of $x$, and define $\textrm{supp}(x)=\{i:x_i\neq 0\}$. The Euclidean norm of $x$ is denoted by $\norm{x}$. Given two functions $\varphi_1:\mathbb{R}_{\ge0}\to\mathbb{R}$ and $\varphi_2:\mathbb{R}_{\ge0}\to\mathbb{R}$, $\varphi_1(n)$ is $O(\varphi_2(n))$ if there exist positive constants $c$ and $N$ such that $|\varphi_1(n)|\le c|\varphi_2(n)|$ for all $n\ge N$. 

\section{Preliminaries}
We begin with some definitions (see, e.g.,\cite{nemhauser1978analysis,conforti1984submodular,bian2017guarantees,kuhnle2018fast}).
\begin{defn}
\label{def:submodularity ratio}
The submodularity ratio of $h:2^{\CS}\to\R_{\ge0}$ is the largest $\gamma\in\R$ such that 
\begin{equation}
\label{eqn:submodularity ratio}
\sum_{v\in \A\setminus \B}\big(h(\{v\}\cup \B)-h(\B)\big)\ge\gamma\big(h(\A\cup\B)-h(\B)\big),
\end{equation}
for all $A,B\subseteq\CS$. The diminishing return (DR) ratio of $h(\cdot)$ is the largest $\kappa\in\R$ such that 
\begin{equation}
\label{eqn:DR ratio}
h(\A\cup\{v\})-h(\A)\ge\kappa\big(h(\B\cup\{v\})-h(\B)\big),
\end{equation}
for all $\A\subseteq\B\subseteq\CS$ and for all $v\in\CS\setminus\B$.
\end{defn}

\begin{defn}
\label{def:curvature} 
The curvature of $h:2^{\CS}\to\R_{\ge0}$ is the smallest $\alpha\in\R$ that satisfies
\begin{equation}
\label{eqn:ex curvature}
h(\A\cup\{v\})-h(\A)\ge(1-\alpha)\big(h(\B\cup\{v\})-h(\B)\big),
\end{equation}
for all $\B\subseteq\A\subseteq\CS$ and for all $v\in\CS\setminus\A$. The extended curvature of $h(\cdot)$ is the smallest $\tilde{\alpha}\in\R$ that satisfies
\begin{equation}
\label{eqn:curvature}
h(\A\cup\{v\})-h(\A)\ge(1-\tilde{\alpha})\big(h(\B\cup\{v\})-h(\B)\big),
\end{equation}
for all $\A,\B\subseteq\CS$ and for all $v\in (\CS\setminus \A)\cap(\CS\setminus \B)$.
\end{defn}
For any monotone nondecreasing set function $h:2^{\CS}\to\R_{\ge0}$, one can check that the submodularity ratio $\gamma$, the DR ratio $\kappa$, the curvature $\alpha$ and the extended curvature $\tilde{\alpha}$ of $h(\cdot)$ satisfy that $\gamma,\kappa,\alpha,\tilde{\alpha}\in[0,1]$. Moreover, we see from Definition~\ref{def:curvature} that $\tilde{\alpha}\ge\alpha$, and it can also been shown that $\gamma\ge\kappa$ (e.g., \cite{bian2017guarantees,kuhnle2018fast,ye2021client}). Further assuming that $h(\cdot)$ is submodular, one can show that $\gamma=\kappa=1$ (e.g., \cite{bian2017guarantees,kuhnle2018fast}). For a modular set function $h:2^{\CS}\to\R_{\ge0}$,\footnote{A set function $h:2^{\CS}\to\R_{\ge0}$ is modular if and only if $h(\A)=\sum_{v\in\A}h(v)$ for all $\A\subseteq\CS$.} we see from Definition~\ref{def:curvature} that the curvature and the extended curvature of $h(\cdot)$ satisfy that $\alpha=\tilde{\alpha}=0$. Thus, the submodularity (resp., DR) ratio of a monotone nondecreasing set function $h(\cdot)$ characterizes the approximate submodularity (resp., approximate DR property) of $h(\cdot)$. The curvatures of $h(\cdot)$ characterize how far the function $h(\cdot)$ is from being modular. Before we proceed, we note that the set $\CS_i$ in \eqref{eqn:obj org} can potentially intersect with $\CS_j$, for any $i,j\in\CS$ with $i\neq j$. Moreover, one can show that a cardinality constraint, a (partitioned) matroid constraint, or multiple budget constraints on the set $\A$ of selected elements are special cases of the constraints in \eqref{eqn:obj org}. In particular, $h_i(\cdot)$ in \eqref{eqn:obj org} reduces to a modular set function for any $i\in[n]$ when considering budget constraints.

\section{Approximation Algorithms}\label{sec:approx alg}
We make the following standing assumption.
\begin{assum}
The set functions $f(\cdot)$ and $h_i(\cdot)$ satisfy that $f(\emptyset)=0$ and $h_i(\emptyset)=0$ for all $i\in[n]$. Further, $h_i(v)>0$ for all $i\in[n]$ and for all $v\in\CS$. 
\end{assum}
Note that \eqref{eqn:obj org} is NP-hard, and cannot be approximated within any constant factor independent of all problem parameters (if P$\neq$NP), even when the constraints in \eqref{eqn:obj org} reduce to a cardinality constraint $|\A|\le H$ \cite{ye2020complexity}. Thus, we aim to provide approximation algorithms for \eqref{eqn:obj org} and characterize the corresponding approximation guarantees in terms of the problem parameters. To simplify the notation in the sequel, for any $\A,\B\subseteq\CS$, we denote
\begin{equation}
\begin{split}
&\delta_{\B}(\A)=f(\A\cup\B)-f(\A)\\
&\delta_{\B}^i(\A) = h_i(\A\cup\B\cap\CS_i)-h_i(\A\cap\CS_i).
\end{split}
\end{equation}
Thus, $\delta_{\B}(\A)$ (resp., $\delta_{\B}^i(\A)$) is the marginal return of $f(\cdot)$ (resp., $h_i(\cdot)$) when adding $\B$ to $\A$.

\subsection{Parallel Greedy Algorithm for a Special Case}\label{sec:greedy algorithms}
We rely on the following assumption and introduce a parallel greedy algorithm (Algorithm~\ref{alg:parallel greedy alg}) for \eqref{eqn:obj org}.
\begin{assum}
\label{ass:disjoint sets}
The ground set $\CS=\cup_{i\in[n]}\CS_i$ in \eqref{eqn:obj org} satisfies that $\CS_i\cap \CS_j=\emptyset$ for all $i,j\in[n]$ with $i\neq j$.
\end{assum}

\begin{algorithm}[tb]
\textbf{Input:} $\CS=\cup_{i\in[n]}\CS_i$, $f:2^{\CS}\to\R_{\ge0}$ and $h_i:2^{\CS}\to\R_{\ge0}$ $\forall i\in[n]$, $H_i\in\R_{\ge0}$ $\forall i\in[n]$
\caption{Parallel greedy algorithm}\label{alg:parallel greedy alg}
\begin{algorithmic}[1]
\For{each $i\in[n]$ in parallel}
\State $\W_i\gets \CS_i$, $\A_i^r\gets\emptyset$
\State $\B_i^r\gets\arg\max_{v\in \CS_i} f(v)$
\While{$\W_i\neq\emptyset$}
    \State $v^{\star}\gets\arg\max_{v\in \W_i}\frac{\delta_v(\A_i^r)}{\delta_v^i(\A_i^r)}$
    \If{$h_i(\A_i^r\cup \{v^{\star}\})\le H_i$}
    \State $\A_i^r\gets \A_i^r\cup \{v^{\star}\}$
    \EndIf
    \State $\W_i\gets \W_i\setminus \{v^{\star}\}$     
    \EndWhile
\State $\A_i^r\gets\arg\max_{\A\in\{\A_i^r,\B_i^r\}}f(\A)$
\EndFor
\State $\A^r\gets \cup_{i\in[n]}\A_i^r$
\State \textbf{Return} $\A^r$
\end{algorithmic}
\end{algorithm}

For each $i\in[n]$ in parallel, Algorithm~\ref{alg:parallel greedy alg} first sets $\CS_i$ to be the ground set $\W_i$ for the algorithm, and then iterates over the current elements in $\W_i$ in the while loop. In particular, the algorithm greedily chooses an element $v\in\W_i$ in line~5 that maximizes the ratio between the marginal returns $\delta_v(\A_i^r)$ and $\delta_v^i(\A_i^r)$ for all $v\in\W_i$. The overall greedy solution is given by $\A^r=\cup_{i\in[n]}\A^r_i$. Thus, one may view Algorithm~\ref{alg:parallel greedy alg} as solving the problem $\max_{\A\subseteq\CS_i}f(\A)$ s.t. $h_i(\A)\le H_i$ for each $i\in[n]$ separately in parallel, and then merge the obtained solutions. Note that the overall time complexity of Algorithm~\ref{alg:parallel greedy alg} is $O((\max_{i\in[n]}|\CS_i|)^2)$. To provide a guarantee on the quality of the approximation for the solution returned by Algorithm~\ref{alg:parallel greedy alg}, we start with the following observation, which follows directly from the definition of the algorithm. 
\begin{observation}
\label{obs:greedy choices}
For any $i\in[n]$ in Algorithm~\ref{alg:parallel greedy alg},  denote $\A_i^r=\{q_1,\dots,q_{|\A_i^r|}\}$ and $\A_{i,j}^r=\{q_1,\dots,q_j\}$ for all $j\in[|\A_i^r|]$ with $\A_{i,0}^r=\emptyset$. Then, there exists $l_i\in[|\A_i^r|]$ such that (1) $q_k\in\arg\max_{v\in \W_i}\frac{\delta_v(\A^r_{i,k-1})}{\delta_v^i(\A_{i,k-1}^r)}$ and $h_i(\A_{i,k}^r)\le H_i$ for all $k\in[l_i]$; and (2) $h_i(\A_{i,l_i}^r\cup \{v^{\star}_{l_i+1}\})>H_i$, where $v^{\star}_{l_i+1}\in\arg\max_{v\in \W_i}\frac{\delta_v(\A_{i,l_i}^r)}{\delta_v^i(\A_{i,l_i}^r)}$.
\end{observation}
\begin{defn}
\label{def:greedy submodularity ratio}
The greedy submodularity ratio of $f:2^{\CS}\to\R_{\ge0}$ in \eqref{eqn:obj org} is the largest $\tilde{\gamma}_f\in\R$ that satisfies $f(\B_i^r)\ge\tilde{\gamma}_f\big(f(\A_{i,l_i}^{r}\cup\{v_{l_i+1}^{\star}\})-f(\A_{i,l_i}^{r})\big)$ for all $i\in[n]$, where $v_{l_i+1}^{\star},\A_{i,l_i}^{r}$ are given in Observation~\ref{obs:greedy choices}, and $\B_i^r$ is given by line~3 of Algorithm~\ref{alg:parallel greedy alg}.
\end{defn}
For monotone nondecreasing $f(\cdot)$ in \eqref{eqn:obj org}, one can check that the greedy submodularity ratio of $f(\cdot)$ satisfies  $\tilde{\gamma}_f\in\R_{\ge0}$. Further assuming that $f(\cdot)$ is submodular, one can show that $\tilde{\gamma}_f\in\R_{\ge1}$. 
\begin{thm}
\label{thm:approx guarantee for parallel greedy}
Suppose that Assumption~\ref{ass:disjoint sets} holds. Let $\A^r$ and $\A^{\star}$ be the solution to \eqref{eqn:obj org} returned by Algorithm~\ref{alg:parallel greedy alg} and an optimal solution to \eqref{eqn:obj org}, respectively. Then,
\begin{multline}
\label{eqn:parallel greedy approx guarantee}
f(\A^r)\ge\frac{(1-\alpha_f)\kappa_f\min\{1,\tilde{\gamma}_f\}}{2}\\\times\min_{i\in[n]}(1-e^{-(1-\tilde{\alpha}_i)\gamma_f}) f(\A^{\star}),
\end{multline}
where $\alpha_f,\kappa_f,\gamma_f\in[0,1]$ and $\tilde{\gamma}_f\in\R_{\ge0}$ are the curvature, DR ratio, submodularity ratio and greedy submodularity ratio of $f(\cdot)$, respectively, and $\tilde{\alpha}_i\in[0,1]$ is the extended curvature of $h_i(\cdot)$ for all $i\in[n]$. 
\end{thm}
We briefly explain the ideas for the proof of Theorem~\ref{thm:approx guarantee for parallel greedy}; a detailed proof is included in Appendix~A. Supposing that $h_i(\cdot)$ is modular for any $i\in[n]$, the choice $v^{\star}$ in line~5 of the algorithm reduces to  $v^{\star}\gets\arg\max_{v\in \W_i}\frac{\delta_v(\A_i^r)}{h_i(v)}$, which is an element $v\in\W_i$ that maximizes the marginal return of the objective function $f(\cdot)$ per unit cost incurred by $h_i(\cdot)$ when adding $v$ to the current greedy solution $\A_i^r$. This renders the greedy nature of the choice $v^{\star}$. To leverage this greedy choice property when $h_i(\cdot)$ is not modular, we use the (extended) curvature of $h_i(\cdot)$ (i.e., $\tilde{\alpha}_i$) to measure how close $h_i(\cdot)$ is to being modular. Moreover, we use $\gamma_f,\tilde{\gamma}_f$ to characterize the approximate submodularity of $f(\cdot)$. Note that the multiplicative factor $(1-\alpha_f)\kappa_f$ in \eqref{eqn:parallel greedy approx guarantee} results from merging $\A_i^r$ for all $i\in[n]$ into $\A^r$ in line~10 of the algorithm.

\begin{rem}
\label{remark:seperatable case}
Under the stronger assumption that the objective function $f(\cdot)$ in \eqref{eqn:obj org} can be written as $f(\A)=\sum_{i\in[n]}f(\A\cap \CS_i)$, similar arguments to those in the proof of Theorem~\ref{thm:approx guarantee for parallel greedy} can be used to show that $f(\A^r)\ge\frac{\min\{1,\tilde{\gamma}_f\}}{2}\min_{i\in[n]}(1-e^{-(1-\tilde{\alpha}_i)\gamma_f}) f(\A^{\star})$.
\end{rem}

\subsection{A Greedy Algorithm for the General Case}\label{sec:greedy algorithm general}
We now introduce a greedy algorithm (Algorithm~\ref{alg:greedy alg}) for general instances of \eqref{eqn:obj org}. The algorithm uses the following definition for the feasible set associated with the constraints in \eqref{eqn:obj org}:
\begin{equation*}
\label{eqn:feasible set}
\CF = \big\{\A\subseteq \CS:h_i(\A\cap\CS_i)\le H_i, \forall i\in[n]\big\}.
\end{equation*}
In the absence of Assumption~\ref{ass:disjoint sets}, Algorithm~\ref{alg:greedy alg} lets $\CS=\cup_{i\in[n]}\CS_i$ be the ground set $\W$ in the algorithm. Algorithm~\ref{alg:greedy alg} then iterates over the current elements in $\W$, and greedily chooses $v\in\W$ and $i\in[n]$ in line~3 such that the ratio between the marginal returns $\delta_v(\A^g)$ and $\delta_v^i(\A^g)$ are maximized for all $v\in\W$ and for all $i\in[n]$. The element $v^{\star}$ will be added to $\A^r$ if the constraint $h_i(\A^g\cup\{v^{\star}\})\le H_i$ is not violated for any $i\in[n]$. Note that different from line~5 in Algorithm~\ref{alg:parallel greedy alg}, the maximization in line~3 of Algorithm~\ref{alg:greedy alg} is also taken with respect to $i\in[n]$. This is because we do not consider the set $\CS_i$ and the constraint associated with $h_i(\cdot)$ for each $i\in[n]$ separately in Algorithm~\ref{alg:greedy alg}. One can check that the time complexity of Algorithm~\ref{alg:greedy alg} is $O(n|\CS|^2)$. In order to characterize the approximation guarantee of Algorithm~\ref{alg:greedy alg}, we introduce the following definition.
\begin{defn}
\label{def:greedy ratio}
Let $\A^g=\{q_1,\dots,q_{|\A^g|}\}$ and $\A^{\star}$ be the solution to \eqref{eqn:obj org} returned by Algorithm~\ref{alg:greedy alg} and an optimal solution to \eqref{eqn:obj org}, respectively. Denote $\A^g_j=\{q_1,\dots,q_j\}$ for all $j\in[|\A^g|]$ with $\A^g_0=\emptyset$. For any $j\in\{0,\dots,|\A^g|-1\}$, the $j$th greedy choice ratio of Algorithm~\ref{alg:greedy alg} is the largest $\psi_j\in\R$ that satisfies
\begin{equation}
\label{eqn:greedy choice ratio}
\frac{\delta_{q_{j+1}}(\A_j^g)}{\delta^{i_j}_{q_{j+1}}(\A^g_j)}\ge \psi_j\frac{\delta_v(\A_j^g)}{\delta_v^i(\A_j^g)}, 
\end{equation}
for all $v\in \A^{\star}\setminus \A_j^g$ and for all $i\in[n]$, where $i_j\in[n]$ is the index of the constraint in \eqref{eqn:obj org} that corresponds to $q_{j+1}$ given by line~3 of Algorithm~\ref{alg:greedy alg}.
\end{defn}
Note that $v^{\star}$ chosen in line~3 of Algorithm~\ref{alg:greedy alg} is not added to the greedy solution $\A_g$ if the constraint in line~4 is violated. Thus, the greedy choice ratio $\psi_j$ given in Definition~\ref{def:greedy ratio} is used to characterize the suboptimality of $q_{j+1}\in\A_g$ in terms of the maximization over $v\in\W$ and $i\in[n]$ in line~3 of the algorithm. Since both $f(\cdot)$ and $h_i(\cdot)$ for all $i\in[n]$ are monotone nondecreasing functions, Definition~\ref{def:greedy ratio} implies that $\psi_j\in\R_{\ge0}$ for all $j\in\{0,\dots,|\A^g|-1\}$.\footnote{The proof of Theorem~\ref{thm:approx guarantee for greedy} shows that we only need to consider the case of both the denominators in \eqref{eqn:greedy choice ratio} being positive.} We also note that for any $j\in\{0,\dots,|\A^g|-1\}$, a lower bound on $\psi_j$ may be obtained by considering all $v\in \CS\setminus \A_j^g$ (instead of $v\in \A^{\star}\setminus \A_j^g$) in Definition~\ref{def:greedy ratio}. Such lower bounds on $\psi_j$ for all $j\in\{0,\dots,|\A^g|-1\}$ can be computed in $O(n|\CS|^2)$ time and in parallel to Algorithm~\ref{alg:greedy alg}. The approximation guarantee for Algorithm~\ref{alg:greedy alg} is provided in the following result proven in Appendix~\ref{app:pf:thm:approx guarantee for greedy}.

\begin{algorithm}[tb]
\textbf{Input:} $\CS=\cup_{i\in[n]}\CS_i$, $f:2^{\CS}\to\R_{\ge0}$ and $h_i:2^{\CS}\to\R_{\ge0}$ $\forall i\in[n]$, $H_i\in\R_{\ge0}$ $\forall i\in[n]$
\caption{Greedy algorithm for general instances of~\eqref{eqn:obj org}}\label{alg:greedy alg}
\begin{algorithmic}[1]
\State $\W\gets\CS$, $\A^g\gets\emptyset$
\While{$\W\neq\emptyset$}
    \State $(v^{\star}, i^{\star})\gets\arg\max_{(v\in \W,i\in[n])}\frac{\delta_v(\A^g)}{\delta_v^i(\A^g)}$
    \If{$(\A^g\cup\{v^{\star}\})\in\CF$}
    \State $\A^g\gets \A^g\cup \{v^{\star}\}$
    \EndIf
    \State $\W\gets \W\setminus \{v^{\star}\}$     
    \EndWhile
\State \textbf{Return} $\A^g$
\end{algorithmic}
\end{algorithm}

\begin{thm}
\label{thm:approx guarantee for greedy}
Let $\A^g$ and $\A^{\star}$ be the solution to \eqref{eqn:obj org} returned by Algorithm~\ref{alg:greedy alg} and an optimal solution to \eqref{eqn:obj org}, respectively. Then,
\begin{equation}
\label{eqn:greedy approx guarantee}
\begin{split}
f(\A^g)&\ge\Big(1-(1-\frac{B}{|\A^g|})^{|\A^g|}\Big)f(\A^{\star})\\
&\ge(1-e^{-B})f(\A^{\star})
\end{split}
\end{equation}
with $B\triangleq\frac{(1-\alpha_h)\gamma_f}{\sum_{i\in[n]}H_i}\sum_{j=0}^{|\A^g|-1}\psi_j\delta^{i_j}_{q_{j+1}}(\A_j^g)$, where $\gamma_f\in[0,1]$ is the submodularity ratio of $f(\cdot)$, $\alpha_h\triangleq\min_{i\in[n]}\tilde{\alpha}_i$ with $\tilde{\alpha}_i\in[0,1]$ to be the extended curvature of $h_i(\cdot)$ for all $i\in[n]$, $\psi_j\in\R_{\ge0}$ is the $j$th greedy choice ratio of Algorithm~\ref{alg:greedy alg} for all $j\in\{0,\dots,|\A^g|-1\}$, and $i_j\in[n]$ is the index of the constraint in \eqref{eqn:obj org} that corresponds to $q_{j+1}$ given by line~3 of Algorithm~\ref{alg:greedy alg}.
\end{thm}

Similarly to Theorem~\ref{thm:approx guarantee for parallel greedy}, we leverage the greedy choice property corresponding to line~3 of Algorithm~\ref{alg:greedy alg} and the properties of $f(\cdot)$ and $h_i(\cdot)$. However, since Algorithm~\ref{alg:greedy alg} considers the constraints associated with $h_i(\cdot)$ for all $i\in[n]$ simultaneously, the proof of Theorem~\ref{thm:approx guarantee for greedy} requires more care, and the approximation guarantee in \eqref{eqn:greedy approx guarantee} does not contain the multiplicative factor $(1-\alpha_f)\kappa_f$.

\begin{rem}
\label{remark:greedy choice ratio}
Similarly to Observation~\ref{obs:greedy choices}, there exists the maximum $l\in[|\A^g|]$ such that for any $q_j\in \A_l^g$, $q_j$ does not violate the condition in line~4 of Algorithm~\ref{alg:greedy alg} when adding to the greedy solution $\A^g$. One can then show via Definition~\ref{def:greedy ratio} that $\psi_{j-1}\ge1$ for all $j\in[|\A^g_l|]$. Further assuming that Assumption~\ref{ass:disjoint sets} holds, one can follow the arguments in the proof of Theorem~\ref{thm:approx guarantee for greedy} and show that 
\begin{equation}
\label{eqn:greedy approx gurantee reduced}
\begin{split}
f(\A^g)&\ge\Big(1-(1-\frac{\tilde{B}}{|\A^g_l|})^{|\A^g_l|}\Big)f(\A^{\star})\\
&\ge(1-e^{-\tilde{B}})f(\A^{\star}),
\end{split}
\end{equation}
where $\tilde{B}\triangleq\frac{(1-\alpha_h)\gamma_f}{\sum_{i\in[n]}H_i}\sum_{i\in[n]}h_i(\A^g_l\cap\CS_i)$.
\end{rem}

\subsection{Comparisons to Existing Results}\label{sec:comparisons to existing results}
Theorems~\ref{thm:approx guarantee for parallel greedy} and \ref{thm:approx guarantee for greedy} generalize several existing results in the literature. First, consider instances of \eqref{eqn:obj org} with a single budget constraint, i.e., $h(\A)=\sum_{v\in\A}h(v)\le H$. We see from Definition~\ref{def:curvature} that the extended curvature of $h(\cdot)$ is $1$. It follows from Remark~\ref{remark:seperatable case} that the approximation guarantee of Algorithm~\ref{alg:parallel greedy alg} provided in Theorem~\ref{thm:approx guarantee for parallel greedy} reduces to $f(\A^r)\ge\frac{\min\{1,\tilde{\gamma}_f\}}{2}(1-e^{-\gamma_f})f(\A^{\star})$, which matches with the approximation guarantee of the greedy algorithm provided in \cite{ye2021parameter}. Further assuming that the objective function $f(\cdot)$ in \eqref{eqn:obj org} is submodular, we have from Definitions~\ref{def:submodularity ratio} and \ref{def:greedy submodularity ratio} that $\gamma_f=1$ and $\tilde{\gamma}_f\ge1$, and the approximation guarantee of Algorithm~\ref{alg:parallel greedy alg} further reduces to $f(\A^r)\ge\frac{1}{2}(1-e^{-1})f(\A^{\star})$, which matches with the results in \cite{khuller1999budgeted,leskovec2007cost}.

Second, consider instances of \eqref{eqn:obj org} with a single cardinality constraint $h(\A)=|\A|\le H$. From Definition~\ref{def:curvature}, we obtain that the curvature of $h(\cdot)$ is $1$. It follows from Remark~\ref{remark:greedy choice ratio} that the approximation guarantee of Algorithm~\ref{alg:greedy alg} provided in Theorem~\ref{thm:approx guarantee for greedy} reduces to $f(\A^g)\ge(1-e^{-\gamma_f})f(\A^{\star})$, which matches with the approximation guarantee of the greedy algorithm provided in \cite{das2018approximate}. Further assuming that the objective function $f(\cdot)$ in \eqref{eqn:obj org} is submodular, we see from Definition~\ref{def:submodularity ratio} that the approximation guarantee of Algorithm~\ref{alg:greedy alg} reduces to $f(\A^g)\ge(1-e^{-1})f(\A^{\star})$, which matches with the result in \cite{nemhauser1978analysis}. 

Third, consider instances of \eqref{eqn:obj org} with a partitioned matroid constraint, i.e., $h_i(\A\cap\CS_i)=|\A\cap\CS_i|\le H_i$ $\forall i\in[n]$. Definition~\ref{def:curvature} shows that the curvature of $h_i(\cdot)$ is $1$ for all $i\in[n]$. Using similar arguments to those in Remark~\ref{remark:greedy choice ratio} and the proof of Theorem~\ref{thm:approx guarantee for greedy}, one can show that Algorithm~\ref{alg:greedy alg} yields the following approximation guarantee:
\begin{equation}
\label{eqn:approx reduce to partitioned matroid}
f(\A^g)\ge\Big(1-(1-\frac{\gamma_f}{\sum_{i\in[n]}H_i})^{|\A^g_l|}\Big)f(\A^{\star}).
\end{equation}
Further assuming that $f(\cdot)$ in \eqref{eqn:obj org} is submodular, i.e., $\gamma_f=1$ in \eqref{eqn:approx reduce to partitioned matroid}, one can check that the approximation guarantee in \eqref{eqn:approx reduce to partitioned matroid} matches with the result in \cite{fisher1978analysis}.

\section{Specific Application Settings}\label{sec:applications}
We now discuss some specific applications that can be captured by the general problem formulation in \eqref{eqn:obj org}. For these applications, we bound the parameters given by Definitions~\ref{def:submodularity ratio}-\ref{def:curvature} and evaluate the resulting approximation guarantees provided in Theorems~\ref{thm:approx guarantee for parallel greedy} and \ref{thm:approx guarantee for greedy}.

\subsection{Sensor Selection}\label{sec:sensor selection problem}
Sensor selection problems arise in many different applications, e.g., \cite{krause2008near,joshi2008sensor,chepuri2014sparsity,ye2020complexity}. A typical scenario is that only a subset of all candidate sensors can be used  to estimate the state of a target environment or system. The goal is to select this subset to optimize an estimation performance metric. If the target system is a dynamical system whose state evolves over time, this problem is sometimes called sensor scheduling, in which different sets of sensors can be selected at different time steps (e.g., \cite{jawaid2015submodularity}) with possibly different constraints on the set of sensors selected at different time steps.

As an example, we can consider the Kalman filtering sensor scheduling (or selection) problem (e.g., \cite{jawaid2015submodularity,tzoumas2016sensor,chamon2017mean,ye2020complexity}) for a linear time-varying system 
\begin{equation}
\label{eqn:LTI}
\begin{split}
x_{k+1} &= A_k x_k + w_k\\
y_k &= C_k x_k + v_k,
\end{split}
\end{equation}
where $A_k\in\R^{n\times n}$, $C_k\in\R^{m\times n}$, $x_0\sim\CN(0,\Pi_0)$ with $\Pi_0\succ 0$, $w_k,v_k$ are zero-mean white Gaussian noise processes with $\BE[w_k w_k^{\top}]=W\succ 0$, $\BE[v_k v_k^{\top}]=\text{diag}(\sigma_{1}^2\cdots\sigma_{m}^2)$, for all $k\in\Z_{\ge0}$ with $\sigma_{i}>0$ for all $i\in[m]$, and $x_0$ is independent of $w_k,v_k$ for all $k\in\Z_{\ge0}$. If there are multiple sensors present, we can let each row in $C_k$ correspond to a candidate sensor at time step $k$. Given a target time step $\ell\in\Z_{\ge0}$, we let $\CS=\{(k,i):i\in[m],k\in\{0,\dots,\ell\}\}$ be the ground set that contains all the candidate sensors at different time steps. Thus, we can write $\CS=\cup_{k\in\{0,\dots,\ell\}}\CS_k$ with $\CS_k=\{(k,i):i\in[m]\}$, where $\CS_k$ (with $|\CS_k|=m$) is the set of sensors available at time step $k$.  Now, for any $\A_k\subseteq\CS_k$, let $C_{\A_k}\in\R^{|\A_k|\times n}$ be the measurement matrix corresponding to the sensors in $\A_k$, i.e., $C_{\A_k}$ contains rows from $C_k$ that correspond to the sensors in $\A_k$. We then consider the following set function $g:2^{\CS}\to\R_{\ge0}$:
\begin{equation}
\label{eqn:def of f in KFSS}
g(\A) = \tr\Big(\big(P_{\ell,\A}^{-1}+\sum_{v\in \A_{\ell}}\sigma_{v}^{-2}C_{v}^{\top} C_{v}\big)^{-1}\Big),
\end{equation}
where  $\A=\cup_{k\in\{0,\dots,\ell\}}\A_k\subseteq\CS$ with $\A_k\subseteq\CS_k$, and $P_{\ell,\A}$ is given recursively via
\begin{equation}
\label{eqn:DARE}
P_{k+1,\A} = W + A_k\big(P_{k,\A}^{-1}+\sum_{v\in\A_k}\sigma_{v}^{-2}C_{v}^{\top}C_{v}\big)^{-1}A_k^{\top},
\end{equation}
for $k=\{0,\dots,\ell-1\}$ with $P_{0,\A}=\Pi_0$. For any $\A\subseteq \CS$, $g(\A)$ is the mean square estimation error of the Kalman filter for estimating the system state $x_{\ell}$ based on the measurements (up until time step $\ell$) from the sensors in $\A$ (e.g., \cite{anderson1979optimal}). Thus, the sensor selection problem can be cast in the framework of \eqref{eqn:obj org} as: 
\begin{equation}
\label{eqn:KFSS}
\begin{split}
&\max_{\A\subseteq\CS}\Big\{f_s(\A)\triangleq g(\emptyset)-g(\A)\Big\}\\
&\textrm{s.t.}\ h_k(\A\cap\CS_k)\le H_k,\ \forall k\in\{0,\dots,\ell\},
\end{split}
\end{equation}
where $H_k\in\R_{\ge0}$, and $h_k(\cdot)$ specifies a constraint on the set of sensors scheduled for any time step $k\in\{0,\dots,\ell\}$. By construction, Assumption~\ref{ass:disjoint sets} holds for problem~\eqref{eqn:KFSS}.

For the objective function, we have the following result for $f_s(\cdot)$; the proof can be adapted from \cite{huber2011optimal,zhang2017sensor,kohara2020sensor} and is omitted for conciseness.
\begin{prop}
\label{prop:properties of f_s}
The set function $f_s:2^{\CS}\to\R_{\ge0}$ in \eqref{eqn:KFSS} is monotone nondecreasing with $f_s(\emptyset)=0$. Moreover, both the submodularity ratio and DR ratio of $f_s(\cdot)$ given in Definition~\ref{def:submodularity ratio} are lower bounded by $\underline{\gamma}$, and both the curvature and extended curvature of $f_s(\cdot)$ given in Definition~\ref{def:curvature} are upper bounded by $\overline{\alpha}$, where $\underline{\gamma}=\frac{\lambda_n(P_{\ell,\emptyset}^{-1})}{\lambda_1(P_{\ell,\CS}^{-1}+\sum_{v\in\CS_{\ell}}\sigma_{v}^{-2}C_v^{\top}C_v)}>0$ and $\overline{\alpha}=1-\underline{\gamma}^2<1$ with $P_{\ell,\emptyset}$ and $P_{\ell,\CS}$ given by Eq.~\eqref{eqn:DARE}.
\end{prop} 
\begin{rem}
Apart from the Kalman filtering sensor selection problem described above, the objective functions in many other formulations such as sensor selection for Gaussian processes \cite{krause2008near}, and sensor selection for hypothesis testing \cite{ye2019sensor} have been shown to be submodular or to have positive submodularity ratio. 
\end{rem}
For the constraints in the sensor selection problems (modeled by $h_k(\cdot)$ in \eqref{eqn:KFSS}), popular choices include a cardinality constraint \cite{krause2008near} or a budget constraint \cite{mo2011sensor,Tzoumas2020LQG} on the set of selected sensors. Our framework can consider such choices individually or simultaneously. More importantly, our framework is general enough to include other relevant constraints. As an example, suppose that the sensors transmit their local information to a fusion center via a (shared) communication channel. Since the fusion center needs to receive the sensor information before the system propagates to the next time step, there are constraints on the communication latency associated with the selected sensors. To ease our presentation, let us consider a specific time step $k\in\{0,\dots,\ell\}$ for the system given by \eqref{eqn:LTI}. Let $[m]$ and $\A\subseteq[m]$ be the set of all the candidate sensors at time step $k$ and the set of sensors selected for time step $k$, respectively. Assume that the sensors in $\A$ transmit the local information to the fusion center using the communication channel in a sequential manner; such an assumption is not restrictive as argued in, e.g., \cite{dinh2020federated}. For any $v\in[m]$, we let $t_v\in\R_{\ge0}$ be the transmission latency corresponding to sensor $v$  when using the communication channel, and let $c_v\in\R_{\ge0}$ be the sensing and computation latency corresponding to sensor $v$. We assume that $t_v,c_v$ are given at the beginning of time step $k$ (e.g., \cite{shi2020joint}). The following assumption says that the sensing and computation latency cannot dominate the transmission latency.
\begin{assum}
\label{ass:time parameters}
For any $v\in[m]$, there exists $r_v\in\R_{>0}$ such that $c_v+t_v-c_u\ge r_v$ for all $u\in[m]$ with $c_u\ge c_v$.
\end{assum}

Note that given a set $\A\subseteq[m]$, the total latency (i.e., the computation and transmission latency) depends on the order in which the sensors in $\A$ transmit. Denote an ordering of the elements in $\A$ as $\hat{\A}=\langle a_1,\dots,a_{|\A|}\rangle$. Define $h^s_c:\hat{\CS}\to\R_{\ge0}$ to be a function that maps a sequence of sensors to the corresponding total latency, where $\hat{\CS}$ is the set that contains all possible sequences of sensors chosen from the set $[m]$. We know from \cite{ye2021client} that the total latency corresponding to $\hat{\A}$ can computed as
\begin{equation}
\label{eqn:recursion for tilde h}
h^s_c(\hat{\A}_j)=\begin{cases}h^s_c(\hat{\A}_{j-1})+t_{a_j}\quad \text{if}\ c_{a_j}<h^s_c(\hat{\A}_{j-1}),\\
c_{a_j}+t_{a_j}\quad \text{if}\ c_{a_j}\ge h_c^s(\hat{\A}_{j-1}),
\end{cases}
\end{equation}
where $\hat{\A}_j=\langle a_1,\dots,a_j \rangle$ for all $j\in[|\A|]$, with $\hat{\A}_0=\emptyset$ and $h^s_c(\emptyset)=0$. 

We further define a set function $h_c:2^{[m]}\to\R_{\ge0}$ such that for any $\A\subseteq[m]$,
\begin{equation}
\label{eqn:def of h}
h_c(\A)=h_c^s(\langle a_1,\dots,a_{|\A|}\rangle),
\end{equation}
where $\langle a_1,\dots,a_{|\A|}\rangle$ orders the elements in $\A$ such that $c_{a_1}\le\cdots\le c_{a_{|\A|}}$. We may now enforce $h_c(\A)\le H$, where $H\in\R_{\ge0}$. Thus, for any $\A\subseteq[m]$, we let the sensors in $\A$ transmit the local information in the order given by \eqref{eqn:def of h} and require the corresponding total latency to be no greater than $H$.\footnote{A similar constraint can be enforced for each time step $k\in\{0,\dots,\ell\}$.} The following result justifies the way $h_c(\cdot)$ orders the selected sensors, and characterizes the curvature of $h_c(\cdot)$.
\begin{prop}
\label{prop:curvature of h_c}
Consider any $\A\subseteq[m]$ and let $\hat{\A}$ be an arbitrary ordering of the elements in $\A$. Then, $h_c(\A)\le h^s_c(\hat{\A})$, where $h(\cdot)$ and $h^s_c(\cdot)$ are defined in \eqref{eqn:recursion for tilde h} and \eqref{eqn:def of h}, respectively. Under Assumption~\ref{ass:time parameters}, it holds that $h_c(\cdot)$ is monotone nondecreasing and that $\tilde{\alpha}_{h_c}\le\tilde{\alpha}_{h_c}^{\prime}\triangleq1-\min_{v\in[m]}\frac{r_v}{t_v}$, where $\tilde{\alpha}_{h_c}\in[0,1]$ is the extended curvature of $h_c(\cdot)$, $r_v\in\R_{>0}$ given in Assumption~\ref{ass:time parameters} and $t_v\in R_{\ge0}$ is the transmission latency corresponding to sensor $v$.
\end{prop}
\begin{pf}
First, for any $\A\subseteq[m]$, denote an arbitrary ordering of the elements in $\A$ as $\hat{\A}=\langle a_1,\dots,a_{|\A|} \rangle$. Then, there exists $\tau\in\Z_{\ge1}$ such that $c_{a_1}\le\cdots\le c_{a_{\tau}}$ and $c_{a_{\tau+1}}\ge c_{a_{\tau}}$, where $c_v\in\R_{\ge0}$ is the computation latency corresponding to sensor $v\in[m]$. Moreover, using the definition of $h^s_c(\cdot)$ in \eqref{eqn:recursion for tilde h}, one can show that 
\begin{equation}
\label{eqn:alternate tilde h}
h_c^s(\hat{\A})=V_c(\hat{\A})+\sum_{j=1}^{|\A|}t_j,
\end{equation}
where $V_c(\hat{\A})\in\R_{\ge0}$ is a function of $\hat{\A}$ that characterizes the time during which the communication channel (shared by all the sensors in $\A$) is idle. Switching the order of $a_{\tau}$ and $a_{\tau+1}$, one can further show that $V_c(\langle a_1,\dots,a_{\tau+1},a_{\tau}\rangle)\le V_c(\langle a_1,\dots,a_{\tau},a_{\tau+1} \rangle)$, which implies via Eq.~\eqref{eqn:alternate tilde h} that $h^s_c(\langle a_1,\dots,a_{\tau+1},a_{\tau}\rangle)\le h^s_c(\langle a_1,\dots,a_{\tau},a_{\tau+1} \rangle)$. It then follows from \eqref{eqn:recursion for tilde h} that $h^s_c(\langle a_1,\dots,a_{\tau-1},a_{\tau+1},a_{\tau},\dots,a_{|\A|} \rangle)\le h^s_c(\hat{\A})$. Repeating the above arguments yields $h_c(\A)\le h^s_c(\hat{\A})$ for any $\A\subseteq[m]$ and any ordering $\hat{\A}$ of the elements in $\A$.

Next, suppose that Assumption~\ref{ass:time parameters} holds. We will show that $h_c(\cdot)$ is monotone nondecreasing and characterize the curvature of $h_c(\cdot)$. To this end, we leverage the expression of $h^s_c(\cdot)$ given by Eq.~\eqref{eqn:alternate tilde h}. Specifically, consider any $\A\subseteq[m]$ and let the elements in $\A$ be ordered such that $\hat{\A}=\langle a_1,\dots,a_{|\A|} \rangle$ with $c_{a_1}\le\cdots\le c_{a_{|\A|}}$. One can first show that $V_c(\hat{\A})=0$. Further considering any $v\in[m]\setminus \A$, one can then show that
\begin{equation*}
\label{eqn:h A union v}
\begin{split}
h_c(\A\cup\{v\})=
\begin{cases}
h_c(\A)+t_v\quad \text{if}\ c_v\ge c_{a_1},\\
h_c(\A)-c_{a_1}+c_v+t_v\quad \text{if}\ c_v < c_{a_1}.
\end{cases}
\end{split}
\end{equation*}
It follows from Assumption~\ref{ass:time parameters} that $h_c(\cdot)$ is monotone nondecreasing, and that $h_c(\A) + r_v\le h_c(\A\cup\{v\})\le h_c(\A) + t_v$, for any $\A\subseteq[m]$ and any $v\in[m]\setminus\A$. Recalling Definition~\ref{def:curvature} completes the proof of the proposition.\hfill$\blacksquare$
\end{pf}

Recalling \eqref{eqn:parallel greedy approx guarantee} (resp., \eqref{eqn:greedy approx gurantee reduced}), substituting $\gamma_f,\kappa_f$ with $\underline{\gamma}$, substituting $\alpha_f$ with $\overline{\alpha}$ from Proposition~\ref{prop:properties of f_s}, and substituting $\tilde{\alpha}_i$ with $\tilde{\alpha}_{h_c}^{\prime}$ from Proposition~\ref{prop:curvature of h_c}, one can obtain the approximation guarantee of Algorithm~\ref{alg:parallel greedy alg} (resp., Algorithm~\ref{alg:greedy alg}) when applied to solve \eqref{eqn:KFSS}.

\begin{rem}
Many other types of constraints in the sensor selection problem can be captured by \eqref{eqn:obj org}. One example is if the selected sensors satisfy certain spatial constraints, e.g.,  two selected sensors may need to be within a certain distance (e.g., \cite{gupta2006stochastic}); or if a mobile robot collects the measurements from the selected sensors (e.g., \cite{prasad2022policies}), the length of the tour of the mobile robot is constrained. Another example is that the budget constraint, where the total cost of the selected sensors is not the sum of the costs of the sensors due to the cost of a sensor can (inversely) depend on the total number of selected sensors (e.g., \cite{iyer2012algorithms}). Such constraints lead to set function constraints on the selected sensors.   
\end{rem}

\subsection{Client Selection for Distributed Optimization}
In a typical distributed optimization framework such as Federated Learning (FL), there is a (central) aggregator and a number of edge devices (i.e., clients) (e.g., \cite{li2019convergence}). Specifically, let $[m]$ be the set of all the candidate clients. For any $v\in[m]$, we assume that the local objective function of client $v$ is given by
\begin{equation*}
\label{eqn:local obj in FL}
F_v(\bw)\triangleq \frac{1}{|\D_v|}\sum_{j=1}^{|\D_v|}\ell_j(\bw;x_{v,j},y_{v,j}),
\end{equation*}
where $\D_v=\{(x_{v,j},y_{v,j}):j\in[|\D_v|]\}$ is the local data set at client $v$, $\ell_j(\cdot)$ is a loss function, and $\bw$ is a model parameter. Here, we let $x_{v,j}\in\R^n$ and $y_{v,j}\in\R$ for all $j\in\D_v$, $\bw\in\R^n$, and $F_v(\bw)\in\R_{\ge0}$ for all $\bw\in\R^n$. The goal is to solve the following global optimization in a distributed manner:
\begin{equation}
\label{eqn:global obj}
\min_{\bw}\Big\{F(\bw)\triangleq\sum_{v\in[m]}\frac{|\D_v|}{D}F_v(\bw)\Big\},
\end{equation}
where $D=\sum_{v\in[m]}|\D_v|$. In general, the FL setup contains multiple rounds of communication between the clients and the aggregator, and solves \eqref{eqn:global obj} using an iterative method (e.g., \cite{li2019convergence}). Specifically, in each round of FL, the aggregator first broadcasts the current global model parameter to the clients. Each client then performs local computations in parallel, in order to update the model parameter using its local dataset via some gradient-based method. Finally, the clients transmit their updated model parameters to the aggregator for global update (see, e.g. \cite{li2019convergence}, for more details).

Similarly to our discussions in Section~\ref{sec:sensor selection problem}, one has to consider constraints (e.g., communication constraints) in FL, which leads to partial participation of the clients (e.g., \cite{reisizadeh2020fedpaq}). Specifically, given a set $\A\subseteq[m]$ of clients that participate in the FL task, it has been shown (e.g., \cite{li2019convergence}) that under certain assumptions on $F(\cdot)$ defined in Eq.~\eqref{eqn:global obj}, the FL algorithm (based on the clients in $\A$) converges to an optimal solution, denoted as $\bw^{\star}_{\A}$, to $\min_{\bw}\Big\{F_{\A}(\bw)\triangleq\sum_{v\in\A}\frac{|\D_v|}{D}F_v(\bw)\Big\}$. We then consider the following client selection problem:
\begin{equation}
\label{eqn:CS for FL}
\begin{split}
&\max_{\A\subseteq\CS}\Big\{f_c(\A)\triangleq F(\bw_{\emptyset})-F(\bw^{\star}_{\A})\Big\}\\
&s.t.\quad h_F(\A)\le T,
\end{split}
\end{equation}
where $\bw_{\emptyset}$ is the initialization of the model parameter.\footnote{We assume that the FL algorithm converges exactly to $w^{\star}_{\A}$ for any $\A\subseteq[m]$. However, the FL algorithm only finds a solution $\tilde{\bw}^{\star}_{\A}$ such that $|F_{\A}(\tilde{\bw}^{\star}_{\A})-F_{\A}(\bw^{\star}_{\A})|=O(1/T_c)$, where $T_c$ is the number of communication rounds between the aggregator and the clients \cite{li2019convergence}. Nonetheless, one can use the techniques in \cite{ye2021parameter} and extend the results for the greedy algorithms provided in this paper to the setting when there are errors in evaluating the objective function $f(\cdot)$ in \eqref{eqn:obj org}.} Similarly to our discussions in Section~\ref{sec:sensor selection problem}, we use $h_F:2^{\CS}\to\R_{\ge0}$ in \eqref{eqn:CS for FL} to characterize the computation latency (for the clients to perform local updates) and the communication latency (for the clients to transmit their local information to the aggregator) in a single round of the FL algorithm.\footnote{We ignore the latency corresponding to the aggregator, since it is typically more powerful than the clients (e.g., \cite{shi2020joint}).} Moreover, we consider the scenario where the clients communicate with the aggregator via a shared channel in a sequential manner. In particular, we may define $h_F(\cdot)$ similarly to $h_c(\cdot)$ given by Eq.~\eqref{eqn:def of h}, and thus the results in Proposition~\ref{prop:curvature of h_c} shown for $h_c(\cdot)$ also hold for $h_F(\cdot)$. The constraint in \eqref{eqn:CS for FL} then ensures that the FL algorithm completes within a certain time limit, when the number of total communication rounds is fixed. Hence, problem~\eqref{eqn:CS for FL} can now be viewed as an instance of problem~\eqref{eqn:obj org}. We also prove the following result for the objective function $f_c(\cdot)$ in problem~\eqref{eqn:CS for FL}.
\begin{prop}
\label{prop:property of f_c}
Suppose that for any $\A\subseteq \B\subseteq[m]$, it holds that (1) $F(\bw^{\star}_{\A})\le F(\bw)$ for all $\bw\in\R^n$ with $\text{supp}(\bw)\subseteq\text{supp}(\bw^{\star}_{\A})$; (2) $\text{supp}(\bw_{\A}^{\star})\subseteq\text{supp}(\bw^{\star}_{\B})$; and (3) $F(\bw^{\star}_{\A})\ge F(\bw_{\B}^{\star})$. Moreover, suppose that for any $v\in[m]$, the local objective function $F_v(\cdot)$ of client $v$ is strongly convex and smooth with parameters $\mu\in\R_{>0}$ and $\rho\in\R_{>0}$, respectively. Then, $f_c(\cdot)$ in problem~\eqref{eqn:CS for FL} is monotone nondecreasing, and both the DR ratio and submodularity ratio of $f_c(\cdot)$ given by Definition~\ref{def:submodularity ratio} are lower bounded by $\mu/\rho$. 
\end{prop}
\begin{pf}
First, since $F(\bw^{\star}_{\A})\ge F(\bw^{\star}_{\B})$ for all $\A\subseteq \B\subseteq[m]$, we see from \eqref{eqn:CS for FL} that $f_c(\cdot)$ is monotone nondecreasing. Next, one can show that the global objective function $F(\cdot)$ defined in \eqref{eqn:global obj} is also strongly convex and smooth with parameters $\mu$ and $\rho$, respectively. That is, for any $\bw_1,\bw_2$ in the domain of $F(\cdot)$, $\frac{\mu}{2}\norm{\bw_2-\bw_1}^2\le-F(\bw_2)+F(\bw_1)+\nabla F(\bw_1)^{\top}(\bw_2-\bw_1)\le\frac{\rho}{2}\norm{\bw_2-\bw_1}^2$. One can now adapt the arguments in the proof of \cite[Theorem~1]{elenberg2018restricted} and show that the bounds on the DR ratio and submodularity ratio of $f_c(\cdot)$ hold. Details of the adaption are omitted here in the interest of space.
\hfill$\blacksquare$
\end{pf}

Recalling \eqref{eqn:greedy approx gurantee reduced}, substituting $\gamma_f$ with $\mu/\rho$ from Proposition~\ref{prop:property of f_c}, and substituting $\tilde{\alpha}_i$ with $\tilde{\alpha}_{h_c}^{\prime}$ from Proposition~\ref{prop:curvature of h_c}, one can obtain the approximation guarantee of Algorithm~\ref{alg:greedy alg} when applied to solve \eqref{eqn:CS for FL}.

One can check that a sufficient condition for the assumptions~(1)-(3) made in Proposition~\ref{prop:property of f_c} to hold is that the datasets from different clients in $[m]$ are non-i.i.d. in the sense that different datasets contain data points with different features, i.e., $\text{supp}(x_{u,i})\cap\text{supp}(x_{v,j})=\emptyset$ for any $u,v\in[m]$ (with $u\neq v$), and for any $i\in\D_u$ and any $j\in\D_v$, where $x_{u,i},x_{v,j}\in\R^n$. In this case, an element $\bw_i$ in the model parameter $\bw\in\R^n$ corresponds to one feature $(x_{u,j})_i$ of the data point $x_{u,j}\in\R^n$, and $\text{supp}(\bw_{\A}^{\star})=\cup_{u\in\A,i\in\D_u}\text{supp}(x_{u,i})$ for any $\A\subseteq[m]$ \cite{elenberg2018restricted}. Note that the datasets from different clients are typically assumed to be non-i.i.d. in FL, since the clients may obtain the local datasets from different data sources \cite{mcmahan2017communication,li2019convergence,shi2020joint}. Moreover, strong convexity and smoothness hold for the loss functions in, e.g., (regularized) linear regression and logistic regression \cite{li2019convergence,shi2020joint}. If the assumptions~(1)-(3) in Proposition~\ref{prop:property of f_c} do not hold, one may use a surrogate for $f_c(\cdot)$ in \eqref{eqn:CS for FL} (see \cite{ye2021client} for more details). 

The FL client selection problem has been studied under various different scenarios (e.g., \cite{nishio2019client,shi2020joint,balakrishnan2021diverse}). In \cite{nishio2019client}, the authors studied a similar client selection problem to the one in this paper, but the objective function considered in \cite{nishio2019client} is simply the sum of sizes of the datasets at the selected clients. In \cite{shi2020joint}, the authors studied a joint optimization problem of bandwidth allocation and client selection, under the training time constraints. However, these works do not provide theoretical performance guarantees for the proposed algorithms. In \cite{balakrishnan2021diverse}, the authors considered a client selection problem with a cardinality constraint on the set of selected clients.

\section{Numerical Results}
\begin{figure}[htbp]
\centering
\subfloat[a][]{
\includegraphics[width=0.48\linewidth]{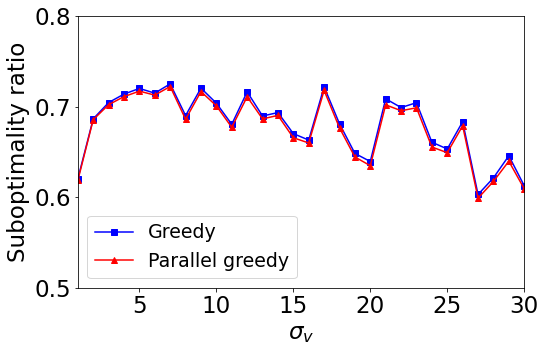}} 
\subfloat[b][]{
\includegraphics[width=0.48\linewidth]{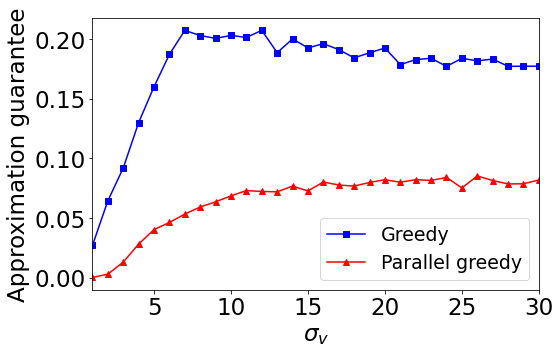}}
\caption{Actual performances and approximation guarantees of Algorithms~\ref{alg:parallel greedy alg}-\ref{alg:greedy alg}.}
\label{fig:real ratio}
\end{figure}

We consider the sensor scheduling problem introduced in \eqref{eqn:KFSS} in Section~\ref{sec:sensor selection problem}, where $h_k(\cdot)$ corresponds to a communication constraint on the set of sensors scheduled for time step $k$ and is defined in Eq.~\eqref{eqn:def of h}, for all $k\in\{0,\dots,\ell\}$. Let the target time step be $\ell=2$, and generate the system matrices $A_k\in\R^{3\times 3}$ and $C_k\in\R^{3\times 3}$ in a random manner, for all $k\in\{0,\dots,\ell\}$. Each row in $C_k$ corresponds to a candidate sensor at time step $k$. We set the input noise covariance as $W=2I_3$, the measurement noise covariance as $\sigma^2_vI_3$ with $\sigma_v\in\{1,\dots,30\}$, and the covariance of $x_0$ as $\Pi_0=I_3$. The ground set $\CS=\{s_{i,k}:i\in[m],k\in\{0,\dots,\ell\}\}$ that contains all the candidate sensors satisfies $|\CS|=3\times 3=9$. For any $k\in\{0,\dots,\ell\}$ and any $i\in[m]$, we generate the computation latency and transmission latency of sensor $s_{i,k}$, denoted as $c_{i,k}$ and $t_{i,k}$, respectively, by sampling exponential distributions with parameters $0.5$ and $0.2$, respectively. Finally, for any $k\in\{0,\dots,\ell\}$, we set the communication constraint to be $H_k=h_k(\CS_k)/2$. We apply Algorithms~\ref{alg:parallel greedy alg} and \ref{alg:greedy alg} to solve the instances of problem~\eqref{eqn:KFSS} constructed above. In Fig.~\ref{fig:real ratio}(a), we plot the actual performances of Algorithms~\ref{alg:parallel greedy alg}-\ref{alg:greedy alg} for $\sigma_v\in\{1,\dots,30\}$, where the actual performance of Algorithm~\ref{alg:parallel greedy alg} (resp., Algorithm~\ref{alg:greedy alg}) is given by $f_s(\A^r)/f_s(\A^{\star})$ (resp., $f_s(\A^g)/f_s(\A^{\star})$), where $f_s(\cdot)$ is given in \eqref{eqn:KFSS}, $\A^r$ (resp., $\A^g$) is the solution to \eqref{eqn:KFSS} returned by Algorithm~\ref{alg:parallel greedy alg} (resp., Algorithm~\ref{alg:greedy alg}), and $\A^{\star}$ is an optimal solution to \eqref{eqn:KFSS} (obtained by brute force). In Fig.~\ref{fig:real ratio}(b), we plot the approximation guarantees of Algorithms~\ref{alg:parallel greedy alg} and \ref{alg:greedy alg} given by Theorems~\ref{thm:approx guarantee for parallel greedy} and \ref{thm:approx guarantee for greedy}, respectively. For any $\sigma_v\in\{1,\dots,30\}$, the results in Fig.~\ref{fig:real ratio}(a)-(b) are averaged over $50$ random instances of problem~\eqref{eqn:KFSS} constructed above. From Fig.~\ref{fig:real ratio}(a), we see that the actual performance of Algorithm~\ref{alg:greedy alg} is slightly better than that of Algorithm~\ref{alg:parallel greedy alg}. We also see that as $\sigma_v$ increases from $1$ to $30$, the actual performances of Algorithms~\ref{alg:parallel greedy alg} and \ref{alg:greedy alg} first tends to be better and then tends to be worse. Compared to Fig.~\ref{fig:real ratio}(a), Fig.~\ref{fig:real ratio}(b) shows that the approximation guarantees of Algorithms~\ref{alg:parallel greedy alg} and \ref{alg:greedy alg} provided by Theorems~\ref{thm:approx guarantee for parallel greedy} and \ref{thm:approx guarantee for greedy}, respectively, are conservative. However, a tighter approximation guarantee potentially yields a better actual performance of the algorithm. 

\begin{figure}[htbp]
\centering
\subfloat[a][]{
\includegraphics[width=0.48\linewidth]{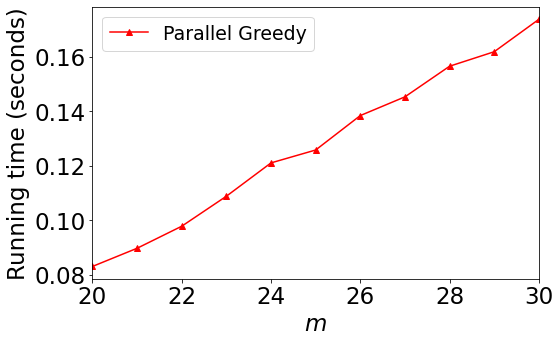}} 
\subfloat[b][]{
\includegraphics[width=0.48\linewidth]{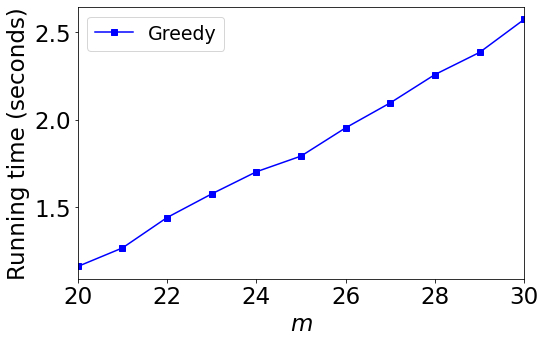}}
\caption{Running times of Algorithms~\ref{alg:parallel greedy alg}-\ref{alg:greedy alg}.}
\label{fig:running times}
\end{figure}

In Fig.~\ref{fig:running times}, we plot the running times of Algorithms~\ref{alg:parallel greedy alg} and \ref{alg:greedy alg} when applied to solve similar random instances of problem~\eqref{eqn:KFSS} to those described above but with $A_k\in\R^{10\times10}$ and $C_k\in\R^{m\times 10}$ for $m\in\{20,\dots30\}$. Note that the simulations are conducted on a Mac with $8$-core CPU, and for any $m\in\{20,\dots,30\}$ the results in Fig.~\ref{fig:running times}(a)-(b) are averaged over $5$ random instances of problem~\eqref{eqn:KFSS}. Fig.~\ref{fig:running times} shows that Algorithm~\ref{alg:parallel greedy alg} runs faster than Algorithm~\ref{alg:greedy alg} matching with our discussions in Section~\ref{sec:approx alg}.

\vspace{-0.3cm}

\section{Conclusion}
We studied the problem of maximizing a monotone nondecreasing set function under multiple constraints, where the constraints are upper bound constraints characterized by monotone nondecreasing set functions. We proposed two greedy algorithms to solve the problem, and analyzed the approximation guarantees of the algorithms, leveraging the notions of submodularity ratio and curvature of set functions. We discussed several important real-world applications of the general problem, and provided bounds on the submodularity ratio and curvature of the set functions in the corresponding instances of the problem. Numerical results show that the two greedy algorithms yield comparable performances that are reasonably good in practice.

\bibliographystyle{plain}        
\bibliography{autosam}           



\appendix
\section{Proof of Theorem~\ref{thm:approx guarantee for parallel greedy}} 
\label{app:pf:thm:approx guarantee for parallel greedy}
\begin{lem}
    Under the setting of Theorem~\ref{thm:approx guarantee for parallel greedy}, we have 
\begin{equation}
\label{eqn:parallel greedy approx guarantee i}
f(\A_i^r)\ge\frac{\min\{1,\tilde{\gamma}_f\}}{2}(1-e^{-(1-\tilde{\alpha}_i)\gamma_f})f(\A_i^{\star}),
\end{equation}
for all $i\in[n]$. 
\end{lem}
\begin{pf}
Under Assumption~\ref{ass:disjoint sets}, we see from the definitions of \eqref{eqn:obj org} and Algorithm~\ref{alg:parallel greedy alg} that $A^r=\cup_{i\in[n]}\A^r_i$ and $\A^{\star}=\cup_{i\in[n]}\A^{\star}_i$, where $\A_i^r,\A_i^{\star}\subseteq \CS_i$ for all $i\in[n]$, and $\A_i^r\cap \A_j^r=\emptyset$, $\A_i^{\star}\cap \A_j^{\star}=\emptyset$ for all $i,j\subseteq[n]$ with $i\neq j$. Now, considering any $i\in[n]$, we note that \eqref{eqn:parallel greedy approx guarantee i} trivially holds if $\tilde{\gamma}_f=0$, $\gamma_f=0$ or $\tilde{\alpha}_i=1$. Thus, in the remaining of this proof, we let $\tilde{\gamma}_f\in\R_{>0}$, $\gamma_f\in(0,1]$ and $\tilde{\alpha}_i\in[0,1)$. Recalling that we have assumed that $h_i(v)>0$ for all $v\in\CS$, we then have from Definition~\ref{def:curvature} that $\delta_v^i(\A)>0$ for all $\A\subseteq\CS$ and for all $v\in \CS\setminus\A$, which implies that line~5 of Algorithm~\ref{alg:parallel greedy alg} is well-defined. Recall from Observation~\ref{obs:greedy choices} that $\A_{i,j}^r=\{q_1,\dots,q_j\}$ for all $j\in[|\A_i^r|]$ with $\A_{i,0}^r=\emptyset$. Moreover, denote $\tilde{\A}_{i,j}=\{\tilde{q}_1,\dots,\tilde{q}_j\}$ for all $j\in[l_i+1]$ with $\tilde{\A}_{i,0}=\emptyset$, where $\tilde{q}_i=q_i$ for all $i\in[l_i]$ and $\tilde{q}_{j+1}=v^{\star}_{l_i+1}$, where $l_i$ and $v_{l_i+1}^{\star}$ are given in Observation~\ref{obs:greedy choices}. Now, considering any $j\in\{0,1,\dots,l_i\}$ and denoting $\tilde{\A}^{\star}_i=\A_i^{\star}\setminus\tilde{\A}_{i,j}^r=\{p_1,\dots,p_{|\tilde{\A}_i^{\star}|}\}$, we have
\begin{align}\nonumber
&f(\A^{\star}_i\cup \tilde{\A}^r_{i,j})-f(\tilde{\A}^r_{i,j})\le\frac{1}{\gamma_f}\sum_{k=1}^{|\tilde{A}^{\star}_i|} \delta_{p_k}(\tilde{\A}_{i,j}^r)\\\nonumber
=&\frac{1}{\gamma_f}\sum_{k=1}^{|\tilde{A}^{\star}_i|}\frac{\delta_{p_k}(\tilde{\A}_{i,j})}{\delta^i_{p_k}(\tilde{\A}_{i,j}^r)}\delta^i_{p_k}(\tilde{\A}_{i,j}^r)\le\frac{\delta_{\tilde{q}_{j+1}}(\tilde{\A}_{i,j}^r)}{\gamma_f\delta^i_{\tilde{q}_{j+1}}(\tilde{\A}_{i,j}^r)}\sum_{k=1}^{|\tilde{\A}^{\star}_i|}\delta^i_{p_k}(\tilde{\A}_{i,j}^r)\\\nonumber
\le&\frac{\delta_{\tilde{q}_{j+1}}(\tilde{A}_{i,j}^r)}{\delta^i_{\tilde{q}_{j+1}}(\tilde{\A}_{i,j}^r)\gamma_f(1-\tilde{\alpha}_i)}\sum_{k=1}^{|\tilde{\A}^{\star}_i|}\delta^i_{p_k}(\{p_1,\dots,p_{k-1}\})\\
=&\frac{\delta_{\tilde{q}_{j+1}}(\tilde{\A}_{i,j}^r)}{(h_i(\tilde{\A}_{i,j+1}^r)-h_i(\tilde{\A}_{i,j}^r))\gamma_f(1-\tilde{\alpha}_i)}h_i(\tilde{\A}_i^{\star}), \label{eqn:alg 1 greedy and opt recursive relation}
\end{align}
where the first inequality follows from Definition~\ref{def:submodularity ratio}, the second inequality follows from the greedy choice in line~5 of Algorithm~\ref{alg:parallel greedy alg}, and the third inequality follows from Definition~\ref{def:curvature}. To proceed, denoting $\Delta_j\triangleq f(\A^{\star}_i)-f(\tilde{\A}^r_{i,j})$ for all $j\in\{0,1,\dots,l_i+1\}$, and noting that $h_i(\tilde{\A}^{\star}_i)\le h_i(\A_i^{\star})\le H_i$, we have from \eqref{eqn:alg 1 greedy and opt recursive relation} the following:
\begin{align}\nonumber
&\Delta_{j+1}\le\big(1-\frac{\gamma_f(1-\tilde{\alpha}_i)(h_i(\tilde{\A}_{i,j+1}^r)-h_i(\tilde{\A}_{i,j}^r))}{H_i}\big)\Delta_j\\\nonumber
\Rightarrow&\Delta_{l_i+1}\le\Delta_0\prod_{j=0}^{l_i} \big(1-\frac{\gamma_f(1-\tilde{\alpha}_i)(h_i(\tilde{A}_{i,j+1}^r)-h_i(\tilde{\A}_{i,j}^r))}{H_i}\big)\\\nonumber
&\qquad\le \Delta_0\prod_{j=0}^{l_i} \big(1-\frac{\gamma_f(1-\tilde{\alpha}_i)(h_i(\tilde{\A}_{i,j+1}^r)-h_i(\tilde{\A}_{i,j}^r))}{h_i(\tilde{\A}_{i,l_i+1}^r)}\big),
\end{align}
where the third inequality follows from $h_i(\tilde{\A}_{i,l_i+1}^r)>H_i$ as we argued above. Note the fact that if $a_1,\dots,a_n\in\R_{>0}$ such that $\sum_{i=1}^n a_i=\alpha G$, where $G\in\R_{>0}$ and $\alpha\in(0,1]$, then the function $\Pi_{i=1}^n(1-\frac{a_i}{G})$ achieves its maximum at $a_1=\cdots=a_n=\frac{\alpha G}{n}$ \cite{kulik2009maximizing}. Since $\sum_{j=0}^{l_i} (h_i(\tilde{\A}_{i,j+1}^r)-h_i(\tilde{\A}_{i,j}^r))=h_i(\tilde{\A}_{i,l_i+1}^r)$ and $\gamma_f(1-\tilde{\alpha}_i)<1$, we have
\begin{align}\nonumber
&\prod_{j=0}^{l_i} \big(1-\frac{\gamma_f(1-\tilde{\alpha}_i)(h_i(\tilde{\A}_{i,j+1}^r)-h_i(\tilde{\A}_{i,j}^r))}{h_i(\tilde{A}_{i,l_i+1}^r)}\big)\\\nonumber
\le&\prod_{j=0}^{l_i} \big(1-\frac{\gamma_f(1-\tilde{\alpha}_i)\frac{h_i(\tilde{\A}_{i,l_i+1}^r)}{l_i+1}}{h_i(\tilde{\A}_{i,l_i+1}^r)}\big)=\big(1-\frac{\gamma_f(1-\tilde{\alpha}_i)}{l_i+1}\big)^{l_i+1}.
\end{align}
It follows that 
\begin{align}\nonumber
&\Delta_{l_i+1}\le \big(1-\frac{\gamma_f(1-\tilde{\alpha}_i)}{l_i+1}\big)^{l_i+1}\Delta_0\le e^{-\gamma_f(1-\tilde{\alpha}_i)}\Delta_0\\
\Rightarrow\ &f(\tilde{\A}_{i,l_i+1}^r)\ge (1-e^{-\gamma_f(1-\tilde{\alpha}_i)})f(\A^{\star}_i).\label{eqn:guarantee for A_l+1}
\end{align}

Recalling from Definition~\ref{def:greedy submodularity ratio} that $f(\tilde{\A}_{i,l_i+1}^r)-f(\tilde{\A}_{i,l_i}^r)\le\frac{1}{\tilde{\gamma}_f}f(\B_i^r)$, we obtain from \eqref{eqn:guarantee for A_l+1} that $f(\tilde{\A}_{i,l_i}^r)+\frac{1}{\tilde{\gamma}_f}f(\B_i^r)\ge(1-e^{-\gamma_f(1-\tilde{\alpha}_i)})f(\A_i^{\star})$. Since $f(\cdot)$ is monotone nondecreasing, it follows  that $f(\A_i^r)+\frac{1}{\tilde{\gamma}_f}f(\B_i^r)\ge(1-e^{-\gamma_f(1-\tilde{\alpha}_i)})f(\A_i^{\star})$, which implies that at least one of $f(\A_i^r)$ and $\frac{1}{\tilde{\gamma}_f}f(\B_i^r)$ is greater than or equal to $\frac{1}{2}(1-e^{-\gamma_f(1-\tilde{\alpha}_i)})f(\A_i^{\star})$. Thus, we see from line~9 in Algorithm~\ref{alg:parallel greedy alg} that \eqref{eqn:parallel greedy approx guarantee i} holds.
\hfill$\blacksquare$
\end{pf} 
{\it Proof of Theorem~\ref{thm:approx guarantee for parallel greedy}:} Since \eqref{eqn:parallel greedy approx guarantee} naturally holds if $\kappa_f=0$, we let $\kappa_f\in(0,1]$ in this proof. Considering $\A^r=\cup_{i\in[n]}\A_i^r$ returned by Algorithm~\ref{alg:parallel greedy alg}, and denoting $\A_1=\cup_{i\in[n-1]}\A_i^r$ and $\A_2=\A_n^r=\{q_1,\dots,q_{|\A_n^r|}\}$, we have
\begin{align}\nonumber
f(\A_1\cup \A_2) &= f(\A_1) + \sum_{j=1}^{|\A_n^r|}\delta_{q_j}(\A_1\cup\{q_1,\dots,q_{j-1}\})\\\nonumber
&\ge f(\A_1) + (1-\alpha_f)\sum_{j=1}^{|\A_n^r|}\delta_{q_j}(\{q_1,\dots,q_{j-1}\})\\
&= f(\A_1) + (1-\alpha_f)f(\A_2).\label{eqn:alg 1 curvature relation}
\end{align}
where the inequality uses the fact that $q_j\notin A_1$ (from Assumption~\ref{ass:disjoint sets}) and Definition~\ref{def:curvature}. Repeating the above arguments for \eqref{eqn:alg 1 curvature relation}, one can show that $f(\cup_{i\in[n]}\A_i^r)\ge(1-\alpha_f)\sum_{i=1}^n f(\A_i^r)$, which implies via \eqref{eqn:parallel greedy approx guarantee i} that
\begin{multline}
\label{eqn:alg 1 curvature and subopt relation}
f(\A^r)\ge\frac{(1-\alpha_f)\min\{1,\tilde{\gamma}_f\}}{2}\\\times\min_{i\in[n]}(1-e^{-(1-\tilde{\alpha}_i)\gamma_i}) \sum_{i=1}^nf(\A_i^{\star}).
\end{multline}
Now, consider the optimal solution $\A^{\star}=\cup_{i\in[n]}\A_i^{\star}$. Using similar arguments to those above, and recalling the definition of the DR ratio $\kappa_f$ of $f(\cdot)$ given in Definition~\ref{def:submodularity ratio}, one can show that $f(\cup_{i\in[n]}\A_i^{\star})\le\frac{1}{\kappa_f}\sum_{i=1}^n f(\A_i^{\star})$,
which together with \eqref{eqn:alg 1 curvature and subopt relation} complete the proof of \eqref{eqn:parallel greedy approx guarantee}.\hfill$\blacksquare$

\section{Proof of Theorem~\ref{thm:approx guarantee for greedy}}
\label{app:pf:thm:approx guarantee for greedy}
First, note that \eqref{eqn:greedy approx guarantee} trivially holds if $\gamma_f=0$, $\alpha_h=1$, or $f(\A^{\star})=0$. Thus, we let $\gamma_f\in(0,1]$, $\tilde{\alpha}_i\in[0,1)$, and $f(\A^{\star})>0$ in this proof. For our analysis in this proof, we assume without loss of generality that $f(\cdot)$ is normalized such that $f(\A^{\star})=1$. Also recalling that we have assumed that $h_i(v)>0$ for all $v\in \CS$, we know from Definition~\ref{def:curvature} that $\delta^i_v(\A)=h_i(\A\cup\{v\})-h_i(\A)>0$ for all $\A\subseteq\CS$ and for all $v\in \CS\setminus\A$, which implies that $\psi_j$ in Definition~\ref{def:greedy ratio} and line~3 of Algorithm~\ref{alg:parallel greedy alg} are well-defined. Denote $\A^g=\{q_1,\dots,q_{|\A^g|}\}$, and $\A^g_j=\{q_1,\dots,q_j\}$ for all $j\in[|\A^g|]$ with $\A^g_0=\emptyset$. We claim that in the remaining of this proof, we can further assume without loss of generality that $\psi_j>0$ for all $j\in\{0,\dots,|\A^g|-1\}$. To prove this claim, we first assume (for contradiction) that there exists $j\in\{0,\dots,|\A^g|-1\}$ such that $\delta_{q_{j+1}}(\A^g_j)=0$. Since $\gamma_f\in(0,1]$ as we argued above, we see from Definition~\ref{def:submodularity ratio} and the greedy choice that $\delta_{\W}(\A^g_j)=0$, where $\W$ is defined and iteratively updated in Algorithm~\ref{alg:greedy alg}. Recalling that $f(\cdot)$ is monotone nondecreasing as we assumed before, we then have that $f(\A^g_j\cup \W_1)=f(\A^g_j\cup \W_2)$ for all $\W_1,\W_2\subseteq \W$. It follows that $\delta_{q_{j^{\prime}+1}}(\A^g_{j^{\prime}})=0$ for all $j^{\prime}\in\{j,\dots,|\A^g|-1\}$. In other words, Algorithm~\ref{alg:greedy alg} is vacuous after adding $q_j$ to the greedy solution $A^g$. Thus, we can assume without loss of generality that $\delta_{q_{j+1}}(\A_j^g)>0$ for all $j\in\{0,\dots,|\A^g|-1\}$, which implies via Definition~\ref{def:greedy ratio} that $\psi_j>0$ for all $j\in\{0,\dots,|\A^g|-1\}$.

To proceed, let us consider any $j\in\{0,\dots,|\A^g|-1\}$. Denoting $\tilde{\A}_j^{\star}=\A^{\star}\setminus \A_j^g$, we have
\begin{align}\nonumber
f(\A^{\star})&\le f(\A_j^g)+\frac{1}{\gamma_f}\sum_{v\in \tilde{\A}^{\star}_j}\delta_v(A^g_j)\\\nonumber
\le f(\A_j^g)&+\frac{1}{\gamma_f}\sum_{i\in[n]}\sum_{v\in \tilde{A}^{\star}_j\cap \CS_i}\frac{\delta_v(\A_j^g)}{\delta_v^i(\A^g_j)}\delta_v^i(\A_j^g)\\
\le f(\A_j^g)&+\frac{\delta_{q_{j+1}}(\A_j^g)}{\gamma_f\psi_j\delta^{i_j}_{q_{j+1}}(\A_j^g)}\sum_{i\in[n]}\sum_{v\in \tilde{\A}^{\star}_j\cap\CS_i}\delta^i_v(\A_j^g),\label{eqn:submodularity and greedy choice ratio}
\end{align}
where the first inequality follows from Definition~\ref{def:submodularity ratio} and the monotonicity of $f(\cdot)$, and the third inequality follows from Definition~\ref{def:greedy ratio}. Denoting $\tilde{\A}^{\star}_j\cap \CS_i=\{p^i_1,\dots,p^i_{|\tilde{\A}^{\star}_j\cap\CS_i|}\}$ for all $i\in[n]$, we further obtain from Definition~\ref{def:curvature} that
\begin{align}\nonumber
&\sum_{i\in[n]}\sum_{v\in \tilde{\A}^{\star}_j\cap\CS_i}\delta^i_v(\A_j^g)\\\nonumber
\le&\frac{1}{1-\alpha_h}\sum_{i\in[n]}\sum_{k=1}^{|\tilde{\A}^{\star}_j\cap S_i|}\delta^i_{p^i_{k}}(\{p_1^i,\dots,p_{k-1}^i\})\\
=&\frac{1}{1-\alpha_h}\sum_{i\in[n]}h_i(\tilde{\A}^{\star}_j\cap\CS_i)\le\frac{\sum_{i\in[n]}H_i}{1-\alpha_h}.\label{eqn:curvature relation}
\end{align}
Denoting 
\begin{equation}
\label{eqn:def of M_j}
M_j = \frac{\sum_{i\in[n]}H_i}{(1-\alpha_h)\gamma_f\psi_j\delta^{i_j}_{q_{j+1}}(\A_j^g)},
\end{equation}
and $\delta_j=\delta_{q_{j+1}}(\A_j^g)$ for all $j\in\{0,\dots,|\A^g|-1\}$, we can combine \eqref{eqn:submodularity and greedy choice ratio}-\eqref{eqn:curvature relation} and obtain that  
\begin{align}
M_j\delta_j + \sum_{k=0}^{j-1}\delta_k\ge 1,\ \forall j\in\{0,\dots,|\A^g|-1\},\label{eqn:lp constraints}
\end{align}
where we use the facts that $f(\A^g_j)=\sum_{k=0}^{j-1}\delta_k$ and $f(\A^{\star})=1$ as we argued before.

In order to prove the approximation guarantee of Algorithm~\ref{alg:greedy alg} given by \eqref{eqn:greedy approx guarantee}, we aim to provide a lower bound on $f(\A^g)/f(\A^{\star})$ and we achieve this by first minimizing $f(\A^g)/f(\A^{\star})=\sum_{j=0}^{|\A^g|-1}\delta_j$ subject to the constraints given in \eqref{eqn:lp constraints}. In other words, we consider the following linear program and its dual:
\begin{equation}
\label{eqn:lp}
\begin{split}
&\min_{\delta_j}\sum_{j=0}^{|\A^g|-1}\delta_j\\
&s.t.\ M_j\delta_j + \sum_{k=0}^{j-1}\delta_k\ge 1,\ \forall j\in\{0,\dots,|\A^g|-1\},
\end{split}
\end{equation}
\begin{equation}
\label{eqn:dual lp}
\begin{split}   
&\max_{\mu_k}\sum_{k=0}^{|\A^g|-1}\mu_k\\
&s.t. \ \mu_k\ge0,\ \forall k\in\{0,\dots,|\A^g|-1\}\\
&\ \ \ M_j\mu_j + \sum_{k=j+1}^{|A^g|-1}\mu_k = 1,\ \forall j\in\{0,\dots,|\A^g|-1\}.
\end{split}
\end{equation}
We will then show that the optimal cost of \eqref{eqn:lp} and \eqref{eqn:dual lp} satisfies \eqref{eqn:greedy approx guarantee}. Noting from Eq.~\eqref{eqn:def of M_j} that $M_j>0$ for all $j\in\{0,\dots,|\A^g|-1\}$, we can obtain the optimal solution to \eqref{eqn:dual lp} by solving the equations in $u_k$ given by the equality constraints in \eqref{eqn:dual lp}, which yields $\mu_k^{\star} = \frac{1}{M_k}\prod_{j=k+1}^{|\A^g|-1}(1-\frac{1}{M_j})$ $\forall k\in\{0,\dots,|\A^g|-1\}$.
Noting from Eq.~\eqref{eqn:def of M_j} that $\sum_{k=0}^{|\A^g|-1}\frac{1}{M_k}=B$, we can further lower bound $\sum_{k=0}^{|\A^g|-1}\mu_k^{\star}$ by solving
\begin{equation}
\label{eqn:lp 2nd}
\begin{split}
\min_{\nu_k}\sum_{k=0}^{|A^g|-1}\nu_k\prod_{j=k+1}^{|\A^g|-1}(1-\nu_j)\quad s.t.\ \sum_{k=0}^{|\A^g|-1}\nu_k = B.
\end{split}
\end{equation}
One can obtain the optimal solution to \eqref{eqn:lp 2nd} by considering a Lagrangian multiplier corresponding to the equality constraint in \eqref{eqn:lp 2nd}, which yields the optimal solution $\nu_k = \frac{B}{|\A^g|}$ $\forall k\in\{0,\dots,|\A^g|-1\}$.
It then follows from our arguments above that 
\begin{align}\nonumber
\sum_{k=0}^{|\A^g|-1}\mu_k^{\star}&\ge\sum_{k=0}^{|\A^g|-1}\frac{B}{|\A^g|}\prod_{j=k+1}^{|\A^g|-1}(1-\frac{B}{|\A^g|})\\\nonumber
&=1-(1-\frac{B}{|\A^g|})^{|\A^g|}\ge1-e^{-B}.\qquad\qquad\quad\hfill\blacksquare
\end{align}
\end{document}